\documentclass[twoside,onecolumn]{article}

\usepackage{blindtext} 

\usepackage[english]{babel} 

\usepackage{titling} 

\usepackage{hyperref} 

\usepackage{amsthm}
\newtheorem{theorem}{Theorem}

\title{Evaluating the Fabius function} 
\author{%
\textsc{Jan Kristian Haugland}\\
\normalsize \href{mailto:admin@neutreeko.net}{admin@neutreeko.net} 
}
\date{\today} 

\renewcommand{\maketitlehookd}{%
	\begin{abstract}
		\noindent The Thue-Morse sequence (1, -1, -1, 1, -1, 1, 1, ...) can in a sense be naturally extended to a continuous function f called the Fabius function. It is shown how to determine the exact value of f(x) whenever x is the ratio between a positive integer and a power of 2.
	\end{abstract}
}

\begin{document}

\maketitle


The Fabius function is a continuous function $f$ on the nonnegative real numbers that satisfies $f(x) + f(1 - x) = 1$ for $0 \leq x \leq 1$, and $$f(x) = \int_0^{2x} f(t) dt$$ for all $x \geq 0$ \cite{Fabius:1966}. Its graph is a bell shape on the interval $\left[ 0, 2 \right]$ that repeats, except for the signs, which follow the Thue-Morse sequence. We will demonstrate how to determine the exact value of $f(x)$ whenever $x$ is the ratio between a positive integer and a power of 2.

It suffices to find closed form expressions for $$2^{2n^2 -2n+1} \left( f \left( {1+x}\over {2^{2n-1}} \right) + f \left( {1-x} \over {2^{2n-1}} \right) \right)$$ and $$2^{2n^2} \left( f \left( {1+x}\over {2^{2n}} \right) - f \left( {1-x} \over {2^{2n}} \right) \right)$$ for $0 \leq x \leq 1$, for each positive integer $n$. Note that for $n=1$, we have $$2 \left( f \left( {1+x}\over {2} \right) + f \left( {1-x} \over {2} \right) \right) = 2$$ and $$4 \left( f \left( {1+x}\over {4} \right) - f \left( {1-x} \over {4} \right) \right) = 4 \left( \int_{0}^{{1+x}\over 2} f(t) dt - \int_{0}^{{1-x}\over 2} f(t) dt \right)$$ $$=4 \int_{{1-x}\over 2}^{{1+x}\over 2} f(t) dt = \int_0^x 2\left( f\left({{1+t}\over 2}\right) + \ f\left( {{1-x}\over 2}\right) \right) dt = 2x$$

\begin{theorem}
	Suppose we have a closed form expression $$a_1 x + a_3 x^3 + a_5 x^5 + ...$$ for $$2^{2n^2} \left( f \left( {1+x}\over {2^{2n}} \right) - f \left( {1-x} \over {2^{2n}} \right) \right)$$ for $0 \leq x \leq 1$, for some positive integer $n$. Then we have $$2^{2n^2 +2n+1}\left( f\left( {{1+x}\over {2^{2n+1}}}\right) + f\left( {{1-x} \over {2^{2n+1}}}\right) \right)$$ $$={1 \over {2^{2n} - 1}} \left( {a_1 \over 3} + {a_3 \over 10} + {a_5 \over 21} + ...\right) + a_1 x^2 + {a_3 \over 2} x^4 + {a_5 \over 3} x^6 + ...$$ and $$2^{2\left( n+1 \right)^2 }\left( f\left( {{1+x}\over {2^{2n+2}}}\right) - f\left( {{1-x} \over {2^{2n+2}}}\right) \right)$$ $$={1 \over {2^{2n} - 1}} \left( {a_1 \over 3} + {a_3 \over 10} + {a_5 \over 21} + ...\right) x + {a_1 \over 3} x^3 + {a_3 \over 10} x^5 + {a_5 \over 21} x^7 + ...$$
\end{theorem}

\begin{proof}
	We have $$2^{2n^2 +2n+1} \left( f \left( {1+x}\over {2^{2n+1}} \right) + f \left( {1-x} \over {2^{2n+1}} \right) \right) $$ $$= 2^{2n^{2} + 2n+1} \left( \int_{0}^{{1+x}\over {2^{2n}}} f(t) dt + \int_{0}^{{1-x}\over {2^{2n}}} f(t) dt \right)$$ $$= 2^{2n^{2} + 2n+1}\left( 2 \int_0^{1 \over {2^{2n}}} f(t) dt + {1 \over {2^{2n}}} \int_0^x \left( f \left( {1+t}\over {2^{2n}} \right) - f \left({1 - t} \over {2^{2n}}\right) \right) \ dt \right)$$ \begin{equation} \label{eq:1} =2^{2n^{2} +2n+2}f \left({1 \over {2^{2n+1}}} \right) + a_1 x^2 + {{a_3}\over 2} x^4 + {a_5 \over 3}x^6 + {a_7 \over 4}x^8 + ... \end{equation} and furthermore $$2^{2n^2 +4n+2} \left( f \left( {1+x}\over {2^{2n+2}} \right) - f \left( {1-x} \over {2^{2n+2}} \right) \right) = 2^{2n^{2} + 4n+2} \int_{{1-x} \over {2^{2n+1}}}^{{1+x}\over {2^{2n+1}}} f(t) dt$$ $$= 2^{2n^{2} + 2n+1} \int_0^{x} \left( f \left( {1 + x} \over {2^{2n+1}} \right) dt + f \left( {1-x}\over {2^{2n+1}} \right) dt \right)$$ \begin{equation} \label{eq:2} =2^{2n^{2} +2n+2}f \left({1 \over {2^{2n+1}}} \right) x + {a_1 \over 3} x^3 + {{a_3}\over 10} x^5 + {a_5 \over 21}x^7 + {a_7 \over 36}x^9 + ... \end{equation}
	On setting $x = 1$ in (2), we can deduce that $$f\left( {1 \over {2^{2n+1}}} \right) = {1 \over {2^{2n^2 +2n+2}}\left( 2^{2n} -1 \right)} \left( {a_1 \over 3}+{a_3 \over 10}+{a_5 \over 21}+{a_7 \over 36}+...\right)$$ and substitute this into (1) and (2), which completes the proof.
\end{proof}

Thus we have, for $0 \leq x \leq 1$: $$2 \left( f \left( {1+x}\over {2} \right) + f \left( {1-x} \over {2} \right) \right) = 2$$ $$2^2 \left( f \left( {1+x}\over {4} \right) - f \left( {1-x} \over {4} \right) \right) = 2x$$ $$2^5 \left( f \left( {1+x}\over {8} \right) + f \left( {1-x} \over {8} \right) \right) = {2 \over 9} + 2x^2$$ $$2^8 \left( f \left( {1+x}\over {16} \right) - f \left( {1-x} \over {16} \right) \right) = {2 \over 9}x + {2 \over 3}x^3$$ $$2^{13} \left( f \left( {1+x}\over {32} \right) + f \left( {1-x} \over {32} \right) \right) = {19 \over 2025} + {2 \over 9} x^2 + {1 \over 3} x^4$$ $$2^{18} \left( f \left( {1+x}\over {64} \right) - f \left( {1-x} \over {64} \right) \right) = {19 \over 2025}x + {2 \over 27}x^3 + {1 \over 15}x^5$$ $$2^{25} \left( f \left( {1+x}\over {128} \right) + f \left( {1-x} \over {128} \right) \right) = {583 \over 2679075} + {19 \over 2025}x^2 + {1 \over 27}x^4 + {1 \over 45}x^6$$ $$2^{32} \left( f \left( {1+x}\over {256} \right) - f \left( {1-x} \over {256} \right) \right) = {583 \over 2679075}x + {19 \over 6075}x^3 + {1 \over 135}x^5 + {1 \over 315}x^7$$ $$2^{41} \left( f \left( {1+x}\over {512} \right) + f \left( {1-x} \over {512} \right) \right)$$ $$= {132809 \over 40989847500} + {583 \over 2679075}x^2 + {19 \over 12150}x^4 + {1 \over 405}x^6 + {1 \over 1260}x^8$$ $$2^{50} \left( f \left( {1+x}\over {1024} \right) - f \left( {1-x} \over {1024} \right) \right)$$ $$= {132809 \over 40989847500}x + {583 \over 8037225}x^3 + {19 \over 60750}x^5 + {1 \over 2835}x^7 + {1 \over 11340}x^9$$ Combining these gives: $$f\left( {1\over 2}\right) = {1 \over 2}$$ $$f\left( {1\over 4}\right) = {5 \over 72}$$ $$f\left( {3\over 4}\right) = {67 \over 72}$$ $$f\left( {1\over 8}\right) = {1 \over 288}$$ $$f\left( {3\over 8}\right) = {73 \over 288}$$ $$f\left( {5\over 8}\right) = {215 \over 288}$$ $$f\left( {7\over 8}\right) = {287 \over 288}$$ $$f\left( {1\over 16}\right) = {143 \over 2073600}$$ $$f\left( {3\over 16}\right) = {46657 \over 2073600}$$ $$f\left( {5\over 16}\right) = {305857 \over 2073600}$$ $$f\left( {7\over 16}\right) = {777743 \over 2073600}$$ $$f\left( {9\over 16}\right) = {1295857 \over 2073600}$$ $$f\left( {11\over 16}\right) = {1767743 \over 2073600}$$ $$f\left( {13\over 16}\right) = {2026943 \over 2073600}$$ $$f\left( {15\over 16}\right) = {2073457 \over 2073600}$$ $$f\left( {1\over 32}\right) = {19 \over 33177600}$$ $$f\left( {3\over 32}\right) = {25219 \over 33177600}$$ $$f\left( {5\over 32}\right) = {334781 \over 33177600}$$ $$f\left( {7\over 32}\right) = {1396781 \over 33177600}$$ $$f\left( {9\over 32}\right) = {3470381 \over 33177600}$$ $$f\left( {11\over 32}\right) = {6555581 \over 33177600}$$ $$f\left( {13\over 32}\right) = {10393219 \over 33177600}$$ $$f\left( {15\over 32}\right) = {14515219 \over 33177600}$$ $$f\left( {17\over 32}\right) = {18662381 \over 33177600}$$ $$f\left( {19\over 32}\right) = {22784381 \over 33177600}$$ $$f\left( {21\over 32}\right) = {26622019 \over 33177600}$$ $$f\left( {23\over 32}\right) = {29707219 \over 33177600}$$ $$f\left( {25\over 32}\right) = {31780819 \over 33177600}$$ $$f\left( {27\over 32}\right) = {32842819 \over 33177600}$$ $$f\left( {29\over 32}\right) = {33152381 \over 33177600}$$ $$f\left( {31\over 32}\right) = {33177581 \over 33177600}$$ and so forth.




\begin{thebibliography}{99} 
	
	\bibitem[Fabius, 1966]{Fabius:1966}
	Fabius, J. (1966).
	\newblock A probabilistic example of a nowhere analytic $C^\infty$-function.
	\newblock {\em Zeitschrift f\"{u}r Wahrscheinlichkeitstheorie und Verwandte Gebiete}, 5:173--174.

\end{thebibliography}
\end{document}